\documentclass[10pt, one side,emlines]{amsart}
\usepackage{tikz}
\usepackage{amssymb,latexsym,xy,eucal,mathrsfs}
\usepackage{bm,amsfonts,amssymb,amsthm,newlfont,eucal, mathrsfs,latexsym, caption}
\textwidth=15.2cm \textheight=24cm \theoremstyle{plain}
\newtheorem{lem}{Lemma}[section]

\newtheorem{thm}[lem]{Theorem}

\newtheorem{cor}[lem]{Corollary}
\newtheorem{nt}[lem]{Note}
\newtheorem{remark}[lem]{Remark}
\theoremstyle{definition}

\newtheorem{obs}[lem]{Observation}

\begin{document}

	\baselineskip 14truept
	\title{Quasi perfect codes in the cartesian product of some graphs}
	
	\author{S. A. Mane and N. V. Shinde }
	
	\dedicatory{Center for Advanced Studies in Mathematics,
		Department of Mathematics,\\ Savitribai Phule Pune University, Pune-411007, India.\\ Department of Mathematics, COEP Technological University, Pune-411005, India.\\
		  manesmruti@yahoo.com : nvs.maths@coeptech.ac.in  \\}
	
	\maketitle
	
	\begin{abstract}

	An important question in the study of quasi-perfect codes is whether such codes can be constructed for all possible lengths $n$. In this paper, we address this question for specific values of $n$. First, we investigate the existence of quasi-perfect codes in the Cartesian product of a graph $G$ and a path (or cycle), assuming that $G$ admits a perfect code. Second, we explore quasi-perfect codes in the Cartesian products of two or three cycles, $C_m\square C_n$ and $C_m\square C_n\square C_l$, as well as in the Cartesian products of two or three paths, $P_m\square P_n$ and $P_m\square P_n\square P_l$.

	\end{abstract}

	
	
	\section{Introduction}

	A graph \( G \) is an ordered pair \( (V(G), E(G)) \), where \( V(G) \) is the set of vertices (or nodes), and \( E(G) \) is the set of edges, each being a two-element subset of \( V(G) \). For any two vertices \( x, y \in V(G) \), the distance \( d(x, y) \) denotes the length of the shortest path between them in the graph.
	
	A \emph{code} \( D \subseteq V(G) \) is a subset of the vertex set. The elements of \( D \) are called \emph{codewords}. A code \( D \) is said to be \emph{\( t \)-error-correcting} if the distance between any two distinct codewords is at least \( 2t + 1 \).
	
	The \emph{covering radius} of a code \( D \) is the smallest integer \( r \) such that for every vertex \( w \in V(G) \), there exists a codeword \( d \in D \) satisfying \( d(w, d) \leq r \). In this case, we say that \( D \) is \emph{\( r \)-covering}.
	
	A code \( D \) is called \emph{\( t \)-perfect} if it is both \( t \)-error-correcting and has covering radius \( t \), i.e., \( r = t \). If the code is perfect with respect to the Lee metric, it is referred to as a \emph{perfect Lee code}. A code is said to be \emph{\( t \)-quasi-perfect} if it is \( t \)-error-correcting and has covering radius \( t + 1 \).
	
	Perfect codes are also known under various terminology in literature: they are referred to as \emph{perfect \( t \)-dominating sets} \cite{li}, \emph{efficient dominating sets} when \( t = 1 \) \cite{li}, and as \emph{perfect distance-\( t \) resource placements} in other contexts \cite{al}, \cite{albda}.

Perfect codes play a central role in the rapidly developing theory of error-correcting codes. However, since perfect codes are relatively rare, the search for and study of \emph{quasi-perfect codes} has emerged as an important area of interest. Both perfect and quasi-perfect codes have been extensively investigated under various metrics, including the Hamming metric~\cite{et}, Lee metric~\cite{le}, and \( \ell_p \)-metric~\cite{zh}.

Quasi-perfect codes with covering radius two or three have been the focus of several studies~\cite{br,dav,davy,et,gab,gor,hor,hora,mo,str,sz,wa,za}. Notably, only a few quasi-perfect codes are known with covering radius greater than three. These include the extended Golay code, which has minimum distance \( 8 \) and covering radius \( 4 \), and the repetition code of length \( 2t \), which has minimum distance \( 2t \) and covering radius \( t \).

The Lee metric was first introduced in~\cite{le} in the context of signal transmission over noisy channels. The study of perfect codes in the Lee metric has been significantly driven by the Golomb–Welch (G–W) conjecture~\cite{go}, which asserts that there exists no perfect Lee \( e \)-error-correcting code of length \( n \) for \( n \geq 3 \) and \( e > 1 \). Although the G–W conjecture remains unresolved, it is widely believed to be true. As a result, research has shifted from seeking perfect Lee codes to exploring \emph{quasi-perfect Lee codes}, which are codes that approximate the properties of perfect ones~\cite{alb}.

Quasi-perfect Lee codes in \( \mathbb{Z}_n \) and \( \mathbb{Z}_{q^n} \) are denoted by \( QPL(n, e) \) and \( QPL(n, e, q) \), respectively. In~\cite{alb}, constructions of \( QPL(2, e, q) \)-codes were presented for all \( e > 1 \) and for all \( q \) satisfying
\[
2e^2 + 2e + 1 \leq q < 2(e + 1)^2 + 2(e + 1) + 1.
\]
In~\cite{ho}, a fast decoding algorithm for these codes was proposed, which operates with constant time complexity. Further, in~\cite{horak}, constructions of \( QPL(n, e) \)-codes for \( n > 2 \) were provided, along with a proof that for any fixed \( n \), there exist only finitely many values of \( e \) for which a linear \( QPL(n, e) \)-code exists. These results suggest that the conditions for the existence of quasi-perfect codes in the Lee metric remain quite restrictive.

In~\cite{zh}, a construction of \( QPL(n, e) \)-codes was presented for (possibly infinitely many) values of \( n \equiv 1 \pmod{6} \). Additionally, quasi-perfect codes were constructed under the \( \ell_p \)-metric.

For further information on quasi-perfect codes, we refer the interested reader to~\cite{al,albd,ba,bi,bu,ca,da,etz,qu,albdai}.

The focus of this paper is the construction of quasi-perfect codes. Given a perfect code in a graph \( G \), we develop a technique to construct a quasi-perfect code in the Cartesian product of \( G \) with certain graphs.

In Section~3, we prove that if \( G \) admits a perfect \( e \)-error-correcting code, then one can construct a quasi-perfect \( e \)-error-correcting code in the Cartesian product \( G \square P_n \) or \( G \square C_n \), where \( P_n \) and \( C_n \) denote the path and cycle graphs, respectively. Furthermore, for \( m, n \geq 3 \) and \( k \geq 1 \), we construct quasi-perfect 2-error-correcting codes in \( P_m \square P_n \square P_{6k - 2} \) and \( C_m \square C_n \square C_{6k} \), using perfect 2-error-correcting codes in \( P_m \square P_n \) and \( C_m \square C_n \), respectively. We also present a construction of a quasi-perfect code in \( P_4 \square P_4 \square P_4 \) based on a perfect code in \( P_2 \square P_2 \square P_2 \). Additionally, we construct quasi-perfect codes in the Cartesian product of two or three cycles, i.e., in \( C_m \square C_n \) and \( C_m \square C_n \square C_l \).

In Section~4, we construct quasi-perfect codes in \( C_n \square C_n \square C_l \) for \( 3 \leq n \leq 19 \) and various values of \( l \), using quasi-perfect codes in \( C_n \square C_n \).

Finally, in Section~5, we present constructions of quasi-perfect codes in the Cartesian products \( P_m \square P_n \) and \( P_m \square P_n \square P_l \), where \( P_m, P_n, P_l \) are path graphs.

	\begin{table}[ht]
		\centering
		\small
		\renewcommand{\arraystretch}{1.2}
		\setlength{\tabcolsep}{6pt}
		\caption{Summary of quasi-perfect code constructions}
		\label{tab:qp-constructions}
		\begin{tabular}{|c|p{3cm}|p{4.2cm}|p{5.3cm}|}
			\hline
			\textbf{Section} & \textbf{Base Graph(s)} & \textbf{Target Graph (Cartesian Product)} & \textbf{Code Properties / Description} \\
			\hline
			3 & \( G \) with a perfect \( e \)-error-correcting code & \( G \square P_n \), \( G \square C_n \) & Construction of quasi-perfect \( e \)-error-correcting codes based on existing perfect codes \\
			\hline
			3 & \( P_m \square P_n \), \( C_m \square C_n \) & \( P_m \square P_n \square P_{6k-2} \), \( C_m \square C_n \square C_{6k} \) & Quasi-perfect 2-error-correcting codes for \( m,n \geq 3,\, k \geq 1 \) \\
			\hline
			3 & \( P_2 \square P_2 \square P_2 \) & \( P_4 \square P_4 \square P_4 \) & Quasi-perfect code constructed using a perfect code in smaller dimension \\
			\hline
			3 & \( C_m, C_n, C_l \) & \( C_m \square C_n \), \( C_m \square C_n \square C_l \) & Quasi-perfect codes in Cartesian products of cycles based on perfect codes \\
			\hline
			4 & \( C_n \square C_n \) & \( C_n \square C_n \square C_l \) & Construction of quasi-perfect codes for \( 3 \leq n \leq 19 \), and various values of \( l \) \\
			\hline
			5 & \( P_m \square P_n \), \( P_m \square P_n \square P_l \) & Same & Quasi-perfect codes in Cartesian products of two and three paths \\
			\hline
		\end{tabular}
	\end{table}

	\section{\textbf{Preliminaries}}

	Throughout this paper, we assume that \( G \) is a simple and connected graph. The symbol \( P_n \) denotes a path on \( n \) vertices. For \( n \geq 3 \), \( C_n \) denotes a cycle on \( n \) vertices; we define \( C_2 \) as a single edge and \( C_1 \) as a single vertex.
	
	For any vertex \( x \in V(G) \), the \emph{ball} of radius \( r \) centered at \( x \) is defined as
	\[
	B_r(x) = \{ y \in V(G) : d(x, y) \leq r \},
	\]
	and the \emph{sphere} of radius \( r \) centered at \( x \) is
	\[
	S_r(x) = \{ y \in V(G) : d(x, y) = r \}.
	\]
	
	For a vertex \( u \in V(G) \) and a subset \( D \subseteq V(G) \), the distance from \( u \) to the set \( D \) is given by
	\[
	d(u, D) = \min\{ d(u, v) : v \in D \}.
	\]
	
	If \( D_1 \) and \( D_2 \) are subsets of vertices from graphs \( G_1 \) and \( G_2 \), respectively, their \emph{direct sum} is defined as
	\[
	D_1 \oplus D_2 = \{ (d_1, d_2) : d_1 \in D_1,\ d_2 \in D_2 \}.
	\]
	
	The \emph{Cartesian product} of two graphs \( G \) and \( H \), denoted \( G \square H \), is the graph with vertex set
	\[
	V(G \square H) = \{ (g, h) : g \in V(G),\ h \in V(H) \},
	\]
	where two vertices \( (g_1, h_1) \) and \( (g_2, h_2) \) are adjacent in \( G \square H \) if and only if either:
	\begin{itemize}
		\item \( h_1 = h_2 \) and \( (g_1, g_2) \in E(G) \), or
		\item \( g_1 = g_2 \) and \( (h_1, h_2) \in E(H) \).
	\end{itemize}
	
	Let \( \mathbb{Z}_q \) denote the ring of integers modulo \( q \), and let \( \mathbb{Z}_q^n \) represent the \( n \)-fold Cartesian product of \( \mathbb{Z}_q \) with itself. For any \( C \subseteq \mathbb{Z}_q^n \) and any \( u \in \mathbb{Z}_q^n \), we define the translate of \( C \) by \( u \) as
	\[
	u + C = \{ u + c : c \in C \}.
	\]
	
	Unless stated otherwise, we assume throughout the paper that all symbols such as \( k, i, j, q, t, n, m \) denote non-negative integers.
	
	For any undefined terminology or notation, we refer the reader to West~\cite{we}.
	
\section {\bf Quasi-perfect codes from perfect codes in Cartesian Products}

	In this section, we construct quasi-perfect \( e \)-error-correcting codes in the Cartesian product \( G \square P_{3k} \) and \( G \square C_{3k} \), using a perfect \( e \)-error-correcting code in the graph \( G \). 
	
	We then investigate the construction of quasi-perfect 2-error-correcting codes in \( P_m \square P_n \square P_{6k-2} \) and \( C_m \square C_n \square C_{6k} \), as well as a quasi-perfect 1-error-correcting code in \( P_4 \square P_4 \square P_4 \). Following this, we explore quasi-perfect code constructions in the Cartesian product of two and three cycles, namely \( C_m \square C_n \) and \( C_m \square C_n \square C_l \).
	
	We begin by constructing quasi-perfect \( e \)-error-correcting codes in the product graph \( G \square P_3 \). Let \( V(G) = \{v_0, v_1, \dots, v_{q-1}\} \); then the vertex set of \( G \square P_n \) can be denoted as
	\[
	V(G \square P_n) = V(G) \oplus \{0, 1, \dots, n - 1\},
	\]
	where \( \oplus \) denotes the direct (Cartesian) product of sets.
	
	\begin{thm}\label{lm1}
		Let \( D \) be a perfect \( e \)-error-correcting code in a graph \( G \). Then:
		
		\begin{itemize}
			\item If \( e = 1 \), there exists a quasi-perfect 1-error-correcting code in the Cartesian product \( G \square P_{3k} \) and \( G \square C_{3k} \) for all \( k \geq 1 \).
			
			\item If \( e \geq 2 \), there exists a quasi-perfect \( e \)-error-correcting code in \( G \square P_3 \) and \( G \square C_3 \).
		\end{itemize}
	\end{thm}
	\begin{proof}
		\begin{itemize}
			\item \textbf{Case 1: \( e \geq 1 \).} Assume that \( D' = D \oplus \{1\} \). Then the minimum distance between any two distinct codewords in \( D' \) remains \( 2e + 1 \), since the product with a singleton set does not affect inter-codeword distances.
			
			Every vertex in \( G \square P_3 \) (or \( G \square C_3 \)) lies within distance at most \( e \) from some codeword in \( D' \), except for vertices in the set \( S_e(x) \oplus \{0,2\} \), for each \( x \in D \). These vertices are at distance \( e + 1 \) from \( D' \). Indeed, for \( x \in D \) and \( u \in S_e(x) \), we have
			\[
			d\big((x,1), (u,0)\big) = d\big((x,1), (u,1)\big) + d\big((u,1), (u,0)\big) = e + 1.
			\]
			Thus, the covering radius of \( D' \) is \( e + 1 \), while the minimum distance is \( 2e + 1 \), meaning \( D' \) is a quasi-perfect \( e \)-error-correcting code in \( G \square P_3 \) (or \( G \square C_3 \)).
			
			\item \textbf{Case 2: \( e = 1 \).} Define
			\[
			D' = \bigcup_{i = 0}^{k - 1} \left(D \oplus \{3i + 1\} \right).
			\]
			From the first part of the proof, each set \( D \oplus \{3i + 1\} \) is a quasi-perfect 1-error-correcting code in the subgraph induced by \( V(G) \oplus \{3i, 3i+1, 3i+2\} \). These subgraphs are disjoint slices of \( G \square P_{3k} \) (or \( G \square C_{3k} \)), and their union covers the entire graph.
			
			Moreover, balls of radius 1 around the codewords in each slice do not overlap with those in other slices. Therefore, \( D' \) is a quasi-perfect 1-error-correcting code in \( G \square P_{3k} \) (or \( G \square C_{3k} \)).
		\end{itemize}
	\end{proof}
	Now, for \( m, n \geq 3 \), we construct quasi-perfect 2-error-correcting codes in the Cartesian product \( P_m \square P_n \square P_{6k-2} \) and \( C_m \square C_n \square C_{6k} \) for all \( k \geq 1 \), using perfect 2-error-correcting codes in \( P_m \square P_n \) and \( C_m \square C_n \), respectively.
	
	\begin{thm}\label{lm2}
		Let \( m, n \geq 3 \). If there exists a perfect 2-error-correcting code in \( P_m \square P_n \), then there exists a quasi-perfect 2-error-correcting code in the Cartesian product \( P_m \square P_n \square P_{6k-2} \) for all \( k \geq 1 \).
	\end{thm}
	
	\begin{proof}
		We begin by constructing a quasi-perfect 2-error-correcting code in \( P_m \square P_n \square P_4 \).
		
		Let \( D_1 \subset V(P_m \square P_n) \) be a perfect 2-error-correcting code in \( P_m \square P_n \). Define a shifted copy \( D_2 = (0,3) + D_1 \). Since \( D_1 \) is perfect, \( D_2 \) also has minimum distance 5, and such a copy always exists due to the structure of the path.
		
		\begin{itemize}
			\item \textbf{Step 1: Construction in \( P_m \square P_n \square P_4 \).} Define
			\[
			D = \left(D_1 \oplus \{0\} \right) \cup \left(D_2 \oplus \{3\} \right).
			\]
			The minimum distance within each layer \( D_1 \oplus \{0\} \) and \( D_2 \oplus \{3\} \) is 5. Also, for any \( u \in D_1 \oplus \{0\} \) and \( v \in D_2 \oplus \{3\} \), we have \( d(u,v) \geq 5 \), because the third coordinate differs by 3, and the base codes are at least distance 5 apart. Thus, the overall minimum distance of \( D \) is 5.
			
			To analyze the covering radius: 
			\begin{itemize}
				\item For each \( w \in D_1 \), the vertices in \( S_2(w) \oplus \{1\} \) are at distance at most 3 from \( D \).
				\item For each \( y \in D_2 \), the vertices in \( S_2(y) \oplus \{2\} \) are also at distance at most 3 from \( D \).
				\item All remaining vertices in \( P_m \square P_n \square P_4 \) are at distance at most 2 from \( D \).
			\end{itemize}
			Therefore, the covering radius of \( D \) is 3, and \( D \) is a quasi-perfect 2-error-correcting code in \( P_m \square P_n \square P_4 \).
			
			\item \textbf{Step 2: Generalization to \( P_m \square P_n \square P_{6k-2} \).} Define
			\[
			C = \bigcup_{i=0}^{k-1} \left( D_1 \oplus \{6i\} \cup D_2 \oplus \{6i + 3\} \right).
			\]
			Each block of length 4 (i.e., from level \( 6i \) to \( 6i+3 \)) replicates the structure of the code \( D \) from Step 1. The blocks are disjoint, and the balls of radius 2 around codewords in each block do not intersect with other blocks. As a result, the minimum distance of \( C \) remains 5, and the covering radius is 3.
			
			Thus, \( C \) is a quasi-perfect 2-error-correcting code in \( P_m \square P_n \square P_{6k-2} \).
		\end{itemize}
	\end{proof}
	
	\medskip
	
	A quasi-perfect 1-error-correcting code in \( P_4 \square P_4 \square P_4 \), constructed from a perfect 1-error-correcting code in \( P_2 \square P_2 \square P_2 \), is illustrated in Figure $1$.
	
	$$\begin{tikzpicture}[scale=0.5]
		
		\draw (0,0) circle [radius=0.1];\draw  (0,5) circle [radius=0.1];\draw  (0,10) circle [radius=0.1];\draw  [fill](0,15) circle [radius=0.2];
		\draw (1,-1) circle [radius=0.1];\draw  (1,4) circle [radius=0.1];\draw  (1,9) circle [radius=0.1];\draw  (1,14) circle [radius=0.1];
		\draw (2,-2) circle [radius=0.1];\draw  (2,3) circle [radius=0.1];\draw (2,8) circle [radius=0.1];\draw  (2,13) circle [radius=0.1];
		\draw  (3,-3) circle [radius=0.1];\draw   [fill](3,2)circle [radius=0.2];\draw  (3,7)circle [radius=0.1];\draw  (3,12) circle [radius=0.1];
		
		\draw [fill](5,0) circle [radius=0.2];\draw  (5,5) circle [radius=0.1];\draw  (5,10) circle [radius=0.1];\draw  (5,15) circle [radius=0.1];
		\draw (6,-1) circle [radius=0.1];\draw  (6,4) circle [radius=0.1];\draw  [fill](6,9) circle [radius=0.2];\draw  (6,14) circle [radius=0.1];
		\draw (7,-2) circle [radius=0.1];\draw  (7,3) circle [radius=0.1];\draw (7,8) circle [radius=0.1];\draw  (7,13) circle [radius=0.1];
		\draw  (8,-3) circle [radius=0.1];\draw   (8,2)circle [radius=0.1];\draw  (8,7)circle [radius=0.1];\draw  (8,12) circle [radius=0.1];
		
		\draw (10,0) circle [radius=0.1];\draw  (10,5) circle [radius=0.1];\draw  (10,10) circle [radius=0.1];\draw  (10,15) circle [radius=0.1];
		\draw (11,-1) circle [radius=0.1];\draw  (11,4) circle [radius=0.1];\draw  (11,9) circle [radius=0.1];\draw  (11,14) circle [radius=0.1];
		\draw (12,-2) circle [radius=0.1];\draw  [fill](12,3) circle [radius=0.2];\draw (12,8) circle [radius=0.1];\draw  (12,13) circle [radius=0.1];
		\draw  (13,-3) circle [radius=0.1];\draw   (13,2)circle [radius=0.1];\draw  (13,7)circle [radius=0.1];\draw  [fill](13,12) circle [radius=0.2];
		
		\draw (15,0) circle [radius=0.1];\draw  (15,5) circle [radius=0.1];\draw  [fill](15,10) circle [radius=0.2];\draw  (15,15) circle [radius=0.1];
		\draw (16,-1) circle [radius=0.1];\draw  (16,4) circle [radius=0.1];\draw  (16,9) circle [radius=0.1];\draw  (16,14) circle [radius=0.1];
		\draw (17,-2) circle [radius=0.1];\draw  (17,3) circle [radius=0.1];\draw (17,8) circle [radius=0.1];\draw  (17,13) circle [radius=0.1];
		\draw [fill] (18,-3) circle [radius=0.2];\draw   (18,2)circle [radius=0.1];\draw  (18,7)circle [radius=0.1];\draw  (18,12) circle [radius=0.1];

		\draw(0,0)--(0,15);\draw(0,0)--(15,0);
		\draw(1,-1)--(1,14);\draw(1,-1)--(16,-1);
		\draw(2,-2)--(2,13);\draw(2,-2)--(17,-2);
		\draw(3,-3)--(3,12);\draw(3,-3)--(18,-3);
		\draw (5,0)--(5,15);\draw (6,-1) --(6,14);\draw (7,-2)--(7,13);\draw  (8,-3)--(8,12);
		\draw (10,0)--(10,15);\draw (11,-1)--(11,14);\draw (12,-2)--(12,13);\draw  (13,-3) -- (13,12) ;
		\draw (15,0)--(15,15);\draw (16,-1)--(16,14);\draw (17,-2)--(17,13);\draw  (18,-3)--(18,12);

		\draw(3,2)--(18,2);\draw(3,7)--(18,7);\draw(3,12)--(18,12);
		\draw(2,3)--(17,3);\draw(2,8)--(17,8);\draw(2,13)--(17,13);
		\draw(1,4)--(16,4);\draw(1,9)--(16,9);\draw(1,14)--(16,14);
		\draw(0,5)--(15,5);\draw(0,10)--(15,10);\draw(0,15)--(15,15);
		
		\draw(0,0)--(3,-3);\draw(5,0)--(8,-3);\draw(10,0)--(13,-3);\draw(15,0)--(18,-3);
		\draw(0,5)--(3,2);\draw(5,5)--(8,2);\draw(10,5)--(13,2);\draw(15,5)--(18,2);
		\draw(0,10)--(3,7);\draw(5,10)--(8,7);\draw(10,10)--(13,7);\draw(15,10)--(18,7);
		\draw(0,15)--(3,12);\draw(5,15)--(8,12);\draw(10,15)--(13,12);\draw(15,15)--(18,12);

		\draw[ultra thick](6,4)--(7,3)--(12,3)--(11,4)--(6,4);
		\draw[ultra thick](6,4)--(6,9);\draw[ultra thick](7,3)--(7,8);\draw[ultra thick](12,3)--(12,8);
		\draw[ultra thick](6,9)--(7,8)--(12,8)--(11,9)--(6,9);
		\draw[ultra thick](11,4)--(11,9);
		
		\draw [thick] (7,-4.5) node[below]{Figure 1. Quasi-perfect $1$-error correcting code in $ P_4\square P_4\square P_4$} ;
		\draw [thick] (7,-5.3) node[below]{(filled circles are codewords)} ;
		
	\end{tikzpicture}$$

	\begin{cor}\label{lm3}
		Let \( m, n \geq 3 \). If there exists a perfect 2-error-correcting code in \( C_m \square C_n \), then for every integer \( k \geq 1 \), there exists a quasi-perfect 2-error-correcting code in the Cartesian product \( C_m \square C_n \square C_{6k} \).
	\end{cor}

\begin{nt}	\label{lm4} \cite{li}
		There exists a perfect 1-error-correcting code in the Cartesian product \( C_3 \square C_6 \square C_2 \), given by the set
		\[
		\{(0,0,0),\, (1,2,0),\, (2,4,0),\, (2,1,1),\, (0,3,1),\, (1,5,1)\}.
		\]
		Using this as a tiling block, one can construct a perfect 1-error-correcting code in \( C_{3p} \square C_{6q} \square C_2 \) for all positive integers \( p \) and \( q \). This tiling approach is illustrated in Figure 3.\end{nt}

\begin{remark}
	In this paper, we present the graphs of the form $C_n \times C_m \times C_l$ (i.e., Cartesian products of cycles) using a series of two-dimensional layers for clarity and ease of understanding. Specifically, we represent the structure as $l$ layers of $C_n \times C_m$, corresponding to the $0$-th through $(l-1)$-th layers along the first coordinate.
	
	To enhance visual clarity and reduce diagrammatic complexity, we intentionally omit the edges that connect the first and last vertices in the $C_n$ and $C_m$ factors. As a result, each layer visually resembles a grid graph $P_n \times P_m$, even though the underlying graph structure remains toroidal (i.e., based on cycles). This omission is purely for illustration purposes and does not affect the correctness or toroidal nature of the graph being represented.
	
	For instance, in Figure~1, we have depicted the graph $P_4 \times P_4 \times P_4$ as a three-dimensional grid structure. The corresponding toroidal graph $C_4 \times C_4 \times C_4$ would appear similar but with additional edges that connect the first and last vertices along each of the three coordinate directions. However, including all such wrap-around edges in the diagram would lead to visual clutter, making it difficult to interpret. Therefore, to maintain readability, we choose to represent the graph using the $P_4 \times P_4 \times P_4$ layout while conceptually referring to the underlying structure as $C_4 \times C_4 \times C_4$ when needed.
\end{remark}

\begin{tikzpicture}[scale=0.8]
	\draw(-1,0) circle [radius=0.1];\draw(-1,1) circle [radius=0.1];\draw(-1,2) circle [radius=0.1];
	\draw[fill](-0.5,0) circle [radius=0.1];\draw(-0.5,1) circle [radius=0.1];\draw(-0.5,2) circle [radius=0.1];
	\draw(0,0) circle [radius=0.1];\draw(0,1) circle [radius=0.1];\draw(0,2) circle [radius=0.1];
	\draw(0.5,0) circle [radius=0.1];\draw(0.5,1) circle [radius=0.1];\draw[fill](0.5,2) circle [radius=0.1];
	\draw(1,0) circle [radius=0.1];\draw(1,1) circle [radius=0.1];\draw(1,2) circle [radius=0.1];
	\draw(1.5,0) circle [radius=0.1];\draw[fill](1.5,1) circle [radius=0.1];\draw(1.5,2) circle [radius=0.1];
	\draw(-1,0)--(-1,2);\draw(-.5,0)--(-.5,2);\draw(0,0)--(0,2);\draw(0.5,0)--(0.5,2);\draw(1,0)--(1,2);\draw(1.5,0)--(1.5,2);
	\draw(-1,0)--(1.5,0);\draw(-1,1)--(1.5,1);\draw(-1,2)--(1.5,2);
	
	\draw(-1,4) circle [radius=0.1];\draw(-1,5) circle [radius=0.1];\draw[fill](-1,6) circle [radius=0.1];
	\draw(-0.5,4) circle [radius=0.1];\draw(-0.5,5) circle [radius=0.1];\draw(-0.5,6) circle [radius=0.1];
	\draw(0,4) circle [radius=0.1];\draw[fill](0,5) circle [radius=0.1];\draw(0,6) circle [radius=0.1];
	\draw(0.5,4) circle [radius=0.1];\draw(0.5,5) circle [radius=0.1];\draw(0.5,6) circle [radius=0.1];
	\draw[fill](1,4) circle [radius=0.1];\draw(1,5) circle [radius=0.1];\draw(1,6) circle [radius=0.1];
	\draw(1.5,4) circle [radius=0.1];\draw(1.5,5) circle [radius=0.1];\draw(1.5,6) circle [radius=0.1];
	\draw(-1,4)--(-1,6);\draw(-.5,4)--(-.5,6);\draw(0,4)--(0,6);\draw(0.5,4)--(0.5,6);\draw(1,4)--(1,6);\draw(1.5,4)--(1.5,6);
	\draw(-1,4)--(1.5,4);\draw(-1,5)--(1.5,5);\draw(-1,6)--(1.5,6);
	\draw [<->] (-2,4) arc (110:250:20pt);
	\draw [<->] (2.5,4) arc (70:-70:20pt);
	\draw [thick] (0,-1.5) node[below]{Figure 2. Perfect $1$-error correcting } ;
	\draw [thick] (0,-2) node[below]{code in $ C_3\times C_6\times C_2$} ;
	\draw [thick] (0,-2.5) node[below]{(filled circles are codewords)} ;
	
	\draw(6,0) circle [radius=0.1];\draw(6,1) circle [radius=0.1];\draw(6,2) circle [radius=0.1];\draw(6,3) circle [radius=0.1];\draw(6,4) circle [radius=0.1];\draw(6,5) circle [radius=0.1];
	\draw[fill](6.5,0) circle [radius=0.1];\draw(6.5,1) circle [radius=0.1];\draw(6.5,2) circle [radius=0.1];\draw[fill](6.5,3) circle [radius=0.1];\draw(6.5,4) circle [radius=0.1];\draw(6.5,5) circle [radius=0.1];
	\draw(7,0) circle [radius=0.1];\draw(7,1) circle [radius=0.1];\draw(7,2) circle [radius=0.1];\draw(7,3) circle [radius=0.1];\draw(7,4) circle [radius=0.1];\draw(7,5) circle [radius=0.1];
	\draw(7.5,0) circle [radius=0.1];\draw(7.5,1) circle [radius=0.1];\draw[fill](7.5,2) circle [radius=0.1];\draw(7.5,3) circle [radius=0.1];\draw(7.5,4) circle [radius=0.1];\draw[fill](7.5,5) circle [radius=0.1];
	\draw(8,0) circle [radius=0.1];\draw(8,1) circle [radius=0.1];\draw(8,2) circle [radius=0.1];\draw(8,3) circle [radius=0.1];\draw(8,4) circle [radius=0.1];\draw(8,5) circle [radius=0.1];
	\draw(8.5,0) circle [radius=0.1];\draw[fill](8.5,1) circle [radius=0.1];\draw(8.5,2) circle [radius=0.1];\draw(8.5,3) circle [radius=0.1];\draw[fill](8.5,4) circle [radius=0.1];\draw(8.5,5) circle [radius=0.1];
	\draw(9,0) circle [radius=0.1];\draw(9,1) circle [radius=0.1];\draw(9,2) circle [radius=0.1];\draw(9,3) circle [radius=0.1];\draw(9,4) circle [radius=0.1];\draw(9,5) circle [radius=0.1];
	\draw[fill](9.5,0) circle [radius=0.1];\draw(9.5,1) circle [radius=0.1];\draw(9.5,2) circle [radius=0.1];\draw[fill](9.5,3) circle [radius=0.1];\draw(9.5,4) circle [radius=0.1];\draw(9.5,5) circle [radius=0.1];
	\draw(10,0) circle [radius=0.1];\draw(10,1) circle [radius=0.1];\draw(10,2) circle [radius=0.1];\draw(10,3) circle [radius=0.1];\draw(10,4) circle [radius=0.1];\draw(10,5) circle [radius=0.1];
	\draw(10.5,0) circle [radius=0.1];\draw(10.5,1) circle [radius=0.1];\draw[fill](10.5,2) circle [radius=0.1];\draw(10.5,3) circle [radius=0.1];\draw(10.5,4) circle [radius=0.1];\draw[fill](10.5,5) circle [radius=0.1];
	\draw(11,0) circle [radius=0.1];\draw(11,1) circle [radius=0.1];\draw(11,2) circle [radius=0.1];\draw(11,3) circle [radius=0.1];\draw(11,4) circle [radius=0.1];\draw(11,5) circle [radius=0.1];
	\draw(11.5,0) circle [radius=0.1];\draw[fill](11.5,1) circle [radius=0.1];\draw(11.5,2) circle [radius=0.1];\draw(11.5,3) circle [radius=0.1];\draw[fill](11.5,4) circle [radius=0.1];\draw(11.5,5) circle [radius=0.1];
	
	\draw(6,7) circle [radius=0.1];\draw(6,8) circle [radius=0.1];\draw[fill](6,9) circle [radius=0.1];\draw(6,10) circle [radius=0.1];\draw(6,11) circle [radius=0.1];\draw[fill](6,12) circle [radius=0.1];
	\draw(6.5,7) circle [radius=0.1];\draw(6.5,8) circle [radius=0.1];\draw(6.5,9) circle [radius=0.1];\draw(6.5,10) circle [radius=0.1];\draw(6.5,11) circle [radius=0.1];\draw(6.5,12) circle [radius=0.1];
	\draw(7,7) circle [radius=0.1];\draw[fill](7,8) circle [radius=0.1];\draw(7,9) circle [radius=0.1];\draw(7,10) circle [radius=0.1];\draw[fill](7,11) circle [radius=0.1];\draw(7,12) circle [radius=0.1];
	\draw(7.5,7) circle [radius=0.1];\draw(7.5,8) circle [radius=0.1];\draw(7.5,9) circle [radius=0.1];\draw(7.5,10) circle [radius=0.1];\draw(7.5,11) circle [radius=0.1];\draw(7.5,12) circle [radius=0.1];
	\draw[fill](8,7) circle [radius=0.1];\draw(8,8) circle [radius=0.1];\draw(8,9) circle [radius=0.1];\draw[fill](8,10) circle [radius=0.1];\draw(8,11) circle [radius=0.1];\draw(8,12) circle [radius=0.1];
	\draw(8.5,7) circle [radius=0.1];\draw(8.5,8) circle [radius=0.1];\draw(8.5,9) circle [radius=0.1];\draw(8.5,10) circle [radius=0.1];\draw(8.5,11) circle [radius=0.1];\draw(8.5,12) circle [radius=0.1];
	\draw(9,7) circle [radius=0.1];\draw(9,8) circle [radius=0.1];\draw[fill](9,9) circle [radius=0.1];\draw(9,10) circle [radius=0.1];\draw(9,11) circle [radius=0.1];\draw[fill](9,12) circle [radius=0.1];
	\draw(9.5,7) circle [radius=0.1];\draw(9.5,8) circle [radius=0.1];\draw(9.5,9) circle [radius=0.1];\draw(9.5,10) circle [radius=0.1];\draw(9.5,11) circle [radius=0.1];\draw(9.5,12) circle [radius=0.1];
	\draw(10,7) circle [radius=0.1];\draw[fill](10,8) circle [radius=0.1];\draw(10,9) circle [radius=0.1];\draw(10,10) circle [radius=0.1];\draw[fill](10,11) circle [radius=0.1];\draw(10,12) circle [radius=0.1];
	\draw(10.5,7) circle [radius=0.1];\draw(10.5,8) circle [radius=0.1];\draw(10.5,9) circle [radius=0.1];\draw(10.5,10) circle [radius=0.1];\draw(10.5,11) circle [radius=0.1];\draw(10.5,12) circle [radius=0.1];
	\draw[fill](11,7) circle [radius=0.1];\draw(11,8) circle [radius=0.1];\draw(11,9) circle [radius=0.1];\draw[fill](11,10) circle [radius=0.1];\draw(11,11) circle [radius=0.1];\draw(11,12) circle [radius=0.1];
	\draw(11.5,7) circle [radius=0.1];\draw(11.5,8) circle [radius=0.1];\draw(11.5,9) circle [radius=0.1];\draw(11.5,10) circle [radius=0.1];\draw(11.5,11) circle [radius=0.1];\draw(11.5,12) circle [radius=0.1];
	
	\draw [thick] (8,-1.5) node[below]{Figure 3. Quasi-perfect $2$-error correcting } ;
	\draw [thick] (8,-2) node[below]{code in $ C_6\square C_{12}\square C_2$ } ;
	\draw [thick] (8,-2.5) node[below]{by using tiling block scheme} ;
	\draw [thick] (8,-3) node[below]{(filled circles are codewords)} ;
	\draw(6,0)--(6,5);\draw(6.5,0)--(6.5,5);\draw(7,0)--(7,5);\draw(7.5,0)--(7.5,5);\draw(8,0)--(8,5);\draw(8.5,0)--(8.5,5);\draw(9,0)--(9,5);\draw(9.5,0)--(9.5,5);\draw(10,0)--(10,5);\draw(10.5,0)--(10.5,5);\draw(11,0)--(11,5);\draw(11.5,0)--(11.5,5);
	\draw(6,0)--(11.5,0);\draw(6,1)--(11.5,1);\draw(6,2)--(11.5,2);\draw(6,3)--(11.5,3);\draw(6,4)--(11.5,4);\draw(6,5)--(11.5,5);
	
	\draw(6,7)--(6,12);\draw(6.5,7)--(6.5,12);\draw(7,7)--(7,12);\draw(7.5,7)--(7.5,12);\draw(8,7)--(8,12);\draw(8.5,7)--(8.5,12);\draw(9,7)--(9,12);\draw(9.5,7)--(9.5,12);\draw(10,7)--(10,12);\draw(10.5,7)--(10.5,12);\draw(11,7)--(11,12);\draw(11.5,7)--(11.5,12);
	\draw(6,7)--(11.5,7);\draw(6,8)--(11.5,8);\draw(6,9)--(11.5,9);\draw(6,10)--(11.5,10);\draw(6,11)--(11.5,11);\draw(6,12)--(11.5,12);
	\draw [<->] (5,6.5) arc (110:250:20pt);
	\draw [<->] (12.85,6.5) arc (70:-70:20pt);
\end{tikzpicture}

\begin{tikzpicture}
	\draw(-1,0) circle [radius=0.1];\draw(-1,1) circle [radius=0.1];\draw(-1,2) circle [radius=0.1];\draw[fill](-1,3) circle [radius=0.1];
	\draw(-0.5,0) circle [radius=0.1];\draw(-0.5,1) circle [radius=0.1];\draw(-0.5,2) circle [radius=0.1];\draw(-0.5,3) circle [radius=0.1];
	\draw(0,0) circle [radius=0.1];\draw(0,1) circle [radius=0.1];\draw(0,2) circle [radius=0.1];\draw(0,3) circle [radius=0.1];
	\draw(0.5,0) circle [radius=0.1];\draw[fill](0.5,1) circle [radius=0.1];\draw(0.5,2) circle [radius=0.1];\draw(0.5,3) circle [radius=0.1];
	\draw(1,0) circle [radius=0.1];\draw(1,1) circle [radius=0.1];\draw(1,2) circle [radius=0.1];\draw(1,3) circle [radius=0.1];
	\draw(1.5,0) circle [radius=0.1];\draw(1.5,1) circle [radius=0.1];\draw(1.5,2) circle [radius=0.1];\draw(1.5,3) circle [radius=0.1];
	\draw(-1,0)--(-1,3);\draw(-.5,0)--(-.5,3);\draw(0,0)--(0,3);\draw(0.5,0)--(0.5,3);\draw(1,0)--(1,3);\draw(1.5,0)--(1.5,3);
	\draw(-1,0)--(1.5,0);\draw(-1,1)--(1.5,1);\draw(-1,2)--(1.5,2);\draw(-1,3)--(1.5,3);
	
	\draw [thick] (0,-1.5) node[below]{Figure 4. Perfect $2$-error correcting } ;
	\draw [thick] (0,-2) node[below]{code in $ C_4\square C_6$} ;
	\draw [thick] (0,-2.5) node[below]{(filled circles are codewords)} ;
	
	\draw(6,0) circle [radius=0.1];\draw(6,1) circle [radius=0.1];\draw(6,2) circle [radius=0.1];\draw(6,3) circle [radius=0.1];\draw[fill](6,4) circle [radius=0.1];
	\draw(6.5,0) circle [radius=0.1];\draw(6.5,1) circle [radius=0.1];\draw(6.5,2) circle [radius=0.1];\draw(6.5,3) circle [radius=0.1];\draw(6.5,4) circle [radius=0.1];
	\draw(7,0) circle [radius=0.1];\draw(7,1) circle [radius=0.1];\draw(7,2) circle [radius=0.1];\draw(7,3) circle [radius=0.1];\draw(7,4) circle [radius=0.1];
	\draw(7.5,0) circle [radius=0.1];\draw(7.5,1) circle [radius=0.1];\draw[fill](7.5,2) circle [radius=0.1];\draw(7.5,3) circle [radius=0.1];\draw(7.5,4) circle [radius=0.1];
	\draw(8,0) circle [radius=0.1];\draw(8,1) circle [radius=0.1];\draw(8,2) circle [radius=0.1];\draw(8,3) circle [radius=0.1];\draw(8,4) circle [radius=0.1];
	\draw(8.5,0) circle [radius=0.1];\draw(8.5,1) circle [radius=0.1];\draw(8.5,2) circle [radius=0.1];\draw(8.5,3) circle [radius=0.1];\draw(8.5,4) circle [radius=0.1];
	\draw(9,0) circle [radius=0.1];\draw(9,1) circle [radius=0.1];\draw(9,2) circle [radius=0.1];\draw(9,3) circle [radius=0.1];\draw(9,4) circle [radius=0.1];
	\draw [thick] (8,-1.5) node[below]{Figure 5. Quasi-perfect $2$-error correcting } ;
	\draw [thick] (8,-2) node[below]{code in $ C_5\square C_7$ by using $2$-perfect code in $ C_4\square C_6$} ;
	\draw [thick] (8,-2.5) node[below]{(filled circles are codewords)} ;
	\draw(6,0)--(6,4);\draw(6.5,0)--(6.5,4);\draw(7,0)--(7,4);\draw(7.5,0)--(7.5,4);\draw(8,0)--(8,4);\draw(8.5,0)--(8.5,4);\draw(9,0)--(9,4);
	\draw(6,0)--(9,0);\draw(6,1)--(9,1);\draw(6,2)--(9,2);\draw(6,3)--(9,3);\draw(6,4)--(9,4);
\end{tikzpicture}

	\begin{thm}\label{lm5}
		A quasi-perfect 1-error-correcting code exists in the Cartesian product \( C_3 \square C_6 \square C_{4k} \), and hence in \( C_{3p} \square C_{6q} \square C_{4k} \) for all positive integers \( p, q, k \).
	\end{thm}
	
	\begin{proof}
		Let 
		\[
		D_0 = \{(0,0), (1,2), (2,4)\}, \quad D_1 = \{(2,1), (0,3), (1,5)\}.
		\]
		Then the set
		\[
		D = (D_0 \oplus \{0\}) \cup (D_1 \oplus \{2\})
		\]
		forms a quasi-perfect 1-error-correcting code in \( C_3 \square C_6 \square C_3 \) and \( C_3 \square C_6 \square C_4 \), since all vertices are within distance 2 from \( D \), and the minimum distance between any two codewords is at least 3.
		
		Now, define:
		\[
		C = \bigcup_{i=0}^{k-1} \left( D_0 \oplus \{4i\} \cup D_1 \oplus \{4i+2\} \right).
		\]
		This union covers \( C_3 \square C_6 \square C_{4k} \) such that each 4-layer segment behaves like the base case above. The codewords remain at minimum pairwise distance 3, and every vertex in the product graph is within distance 2 of some codeword. Hence, \( C \) is a quasi-perfect 1-error-correcting code in \( C_3 \square C_6 \square C_{4k} \).
		
		Finally, since \( C_3 \square C_6 \square C_{4k} \) tiles \( C_{3p} \square C_{6q} \square C_{4k} \) for all positive integers \( p, q \), we can extend the code to the larger product by periodic repetition, preserving both minimum distance and covering radius. Thus, a quasi-perfect 1-error-correcting code also exists in \( C_{3p} \square C_{6q} \square C_{4k} \).
	\end{proof}
	
	\medskip
	
	Using the technique from Theorem~\ref{lm5}, we obtain the following generalization:
	
	\begin{thm}\label{lm6}
		Let \( m, n \geq 2 \). If there exists a perfect 1-error-correcting code in the Cartesian product \( C_m \square C_n \square C_2 \), then a quasi-perfect 1-error-correcting code exists in \( C_m \square C_n \square C_4 \), and hence in \( C_m \square C_n \square C_{4k} \) for all positive integers \( k \).
	\end{thm}
\begin{tikzpicture}[scale=0.7]
	\draw(-.5,5) circle [radius=0.1];\draw(-.5,6) circle [radius=0.1];\draw(-.5,7) circle [radius=0.1];\draw[fill](-.5,8) circle [radius=0.1];
	\draw(0,5) circle [radius=0.1];\draw(0,6) circle [radius=0.1];\draw(0,7) circle [radius=0.1];\draw(0,8) circle [radius=0.1];
	\draw(.5,5) circle [radius=0.1];\draw(.5,6) circle [radius=0.1];\draw(.5,7) circle [radius=0.1];\draw(.5,8) circle [radius=0.1];
	\draw(1,5) circle [radius=0.1];\draw(1,6) circle [radius=0.1];\draw(1,7) circle [radius=0.1];\draw(1,8) circle [radius=0.1];
	\draw(1.5,5) circle [radius=0.1];\draw(1.5,6) circle [radius=0.1];\draw(1.5,7) circle [radius=0.1];\draw(1.5,8) circle [radius=0.1];
	\draw(2,5) circle [radius=0.1];\draw(2,6) circle [radius=0.1];\draw(2,7) circle [radius=0.1];\draw(2,8) circle [radius=0.1];
	\draw(2.5,5) circle [radius=0.1];\draw(2.5,6) circle [radius=0.1];\draw(2.5,7) circle [radius=0.1];\draw(2.5,8) circle [radius=0.1];
	\draw(3,5) circle [radius=0.1];\draw(3,6) circle [radius=0.1];\draw(3,7) circle [radius=0.1];\draw(3,8) circle [radius=0.1];
	\draw(3.5,5) circle [radius=0.1];\draw(3.5,6) circle [radius=0.1];\draw(3.5,7) circle [radius=0.1];\draw[fill](3.5,8) circle [radius=0.1];
	\draw(4,5) circle [radius=0.1];\draw(4,6) circle [radius=0.1];\draw(4,7) circle [radius=0.1];\draw(4,8) circle [radius=0.1];
	\draw(4.5,5) circle [radius=0.1];\draw(4.5,6) circle [radius=0.1];\draw(4.5,7) circle [radius=0.1];\draw(4.5,8) circle [radius=0.1];
	\draw(5,5) circle [radius=0.1];\draw(5,6) circle [radius=0.1];\draw(5,7) circle [radius=0.1];\draw(5,8) circle [radius=0.1];
	\draw(5.5,5) circle [radius=0.1];\draw(5.5,6) circle [radius=0.1];\draw(5.5,7) circle [radius=0.1];\draw(5.5,8) circle [radius=0.1];
	\draw(6,5) circle [radius=0.1];\draw(6,6) circle [radius=0.1];\draw(6,7) circle [radius=0.1];\draw(6,8) circle [radius=0.1];
	\draw(6.5,5) circle [radius=0.1];\draw(6.5,6) circle [radius=0.1];\draw(6.5,7) circle [radius=0.1];\draw(6.5,8) circle [radius=0.1];
	\draw(7,5) circle [radius=0.1];\draw(7,6) circle [radius=0.1];\draw(7,7) circle [radius=0.1];\draw(7,8) circle [radius=0.1];
	\draw [thick] (3,4.5) node[below]{$0^{th}$ layer} ;

	\draw(6,0) circle [radius=0.1];\draw(6,1) circle [radius=0.1];\draw(6,2) circle [radius=0.1];\draw(6,3) circle [radius=0.1];
	\draw(6.5,0) circle [radius=0.1];\draw(6.5,1) circle [radius=0.1];\draw(6.5,2) circle [radius=0.1];\draw(6.5,3) circle [radius=0.1];
	\draw(7,0) circle [radius=0.1];\draw(7,1) circle [radius=0.1];\draw(7,2) circle [radius=0.1];\draw(7,3) circle [radius=0.1];
	\draw(7.5,0) circle [radius=0.1];\draw(7.5,1) circle [radius=0.1];\draw(7.5,2) circle [radius=0.1];\draw(7.5,3) circle [radius=0.1];
	\draw(8,0) circle [radius=0.1];\draw(8,1) circle [radius=0.1];\draw(8,2) circle [radius=0.1];\draw(8,3) circle [radius=0.1];
	\draw(8.5,0) circle [radius=0.1];\draw(8.5,1) circle [radius=0.1];\draw(8.5,2) circle [radius=0.1];\draw(8.5,3) circle [radius=0.1];
	\draw(9,0) circle [radius=0.1];\draw(9,1) circle [radius=0.1];\draw(9,2) circle [radius=0.1];\draw(9,3) circle [radius=0.1];
	\draw(9.5,0) circle [radius=0.1];\draw(9.5,1) circle [radius=0.1];\draw(9.5,2) circle [radius=0.1];\draw(9.5,3) circle [radius=0.1];
	\draw(10,0) circle [radius=0.1];\draw(10,1) circle [radius=0.1];\draw(10,2) circle [radius=0.1];\draw(10,3) circle [radius=0.1];
	\draw(10.5,0) circle [radius=0.1];\draw(10.5,1) circle [radius=0.1];\draw(10.5,2) circle [radius=0.1];\draw(10.5,3) circle [radius=0.1];
	\draw(11,0) circle [radius=0.1];\draw(11,1) circle [radius=0.1];\draw(11,2) circle [radius=0.1];\draw(11,3) circle [radius=0.1];
	\draw(11.5,0) circle [radius=0.1];\draw(11.5,1) circle [radius=0.1];\draw(11.5,2) circle [radius=0.1];\draw(11.5,3) circle [radius=0.1];
	\draw(12,0) circle [radius=0.1];\draw(12,1) circle [radius=0.1];\draw(12,2) circle [radius=0.1];\draw(12,3) circle [radius=0.1];
	\draw(12.5,0) circle [radius=0.1];\draw(12.5,1) circle [radius=0.1];\draw(12.5,2) circle [radius=0.1];\draw(12.5,3) circle [radius=0.1];
	\draw(13,0) circle [radius=0.1];\draw(13,1) circle [radius=0.1];\draw(13,2) circle [radius=0.1];\draw(13,3) circle [radius=0.1];
	\draw(13.5,0) circle [radius=0.1];\draw(13.5,1) circle [radius=0.1];\draw(13.5,2) circle [radius=0.1];\draw(13.5,3) circle [radius=0.1];
	\draw [thick] (10,-.5) node[below]{$1^{st}$ layer} ;
	
	\draw(6,10) circle [radius=0.1];\draw(6,11) circle [radius=0.1];\draw(6,12) circle [radius=0.1];\draw(6,13) circle [radius=0.1];
	\draw(6.5,10) circle [radius=0.1];\draw(6.5,11) circle [radius=0.1];\draw(6.5,12) circle [radius=0.1];\draw(6.5,13) circle [radius=0.1];
	\draw(7,10) circle [radius=0.1];\draw(7,11) circle [radius=0.1];\draw(7,12) circle [radius=0.1];\draw(7,13) circle [radius=0.1];
	\draw(7.5,10) circle [radius=0.1];\draw(7.5,11) circle [radius=0.1];\draw(7.5,12) circle [radius=0.1];\draw(7.5,13) circle [radius=0.1];
	\draw(8,10) circle [radius=0.1];\draw(8,11) circle [radius=0.1];\draw(8,12) circle [radius=0.1];\draw(8,13) circle [radius=0.1];
	\draw(8.5,10) circle [radius=0.1];\draw(8.5,11) circle [radius=0.1];\draw(8.5,12) circle [radius=0.1];\draw(8.5,13) circle [radius=0.1];
	\draw(9,10) circle [radius=0.1];\draw(9,11) circle [radius=0.1];\draw(9,12) circle [radius=0.1];\draw(9,13) circle [radius=0.1];
	\draw(9.5,10) circle [radius=0.1];\draw(9.5,11) circle [radius=0.1];\draw(9.5,12) circle [radius=0.1];\draw(9.5,13) circle [radius=0.1];
	\draw(10,10) circle [radius=0.1];\draw(10,11) circle [radius=0.1];\draw(10,12) circle [radius=0.1];\draw(10,13) circle [radius=0.1];
	\draw(10.5,10) circle [radius=0.1];\draw(10.5,11) circle [radius=0.1];\draw(10.5,12) circle [radius=0.1];\draw(10.5,13) circle [radius=0.1];
	\draw(11,10) circle [radius=0.1];\draw(11,11) circle [radius=0.1];\draw(11,12) circle [radius=0.1];\draw(11,13) circle [radius=0.1];
	\draw(11.5,10) circle [radius=0.1];\draw(11.5,11) circle [radius=0.1];\draw(11.5,12) circle [radius=0.1];\draw(11.5,13) circle [radius=0.1];
	\draw(12,10) circle [radius=0.1];\draw(12,11) circle [radius=0.1];\draw(12,12) circle [radius=0.1];\draw(12,13) circle [radius=0.1];
	\draw(12.5,10) circle [radius=0.1];\draw(12.5,11) circle [radius=0.1];\draw(12.5,12) circle [radius=0.1];\draw(12.5,13) circle [radius=0.1];
	\draw(13,10) circle [radius=0.1];\draw(13,11) circle [radius=0.1];\draw(13,12) circle [radius=0.1];\draw(13,13) circle [radius=0.1];
	\draw(13.5,10) circle [radius=0.1];\draw(13.5,11) circle [radius=0.1];\draw(13.5,12) circle [radius=0.1];\draw(13.5,13) circle [radius=0.1];
	\draw [thick] (10,9.5) node[below]{$3^{rd}$ layer} ;
	
	\draw(12.5,5) circle [radius=0.1];\draw(12.5,6) circle [radius=0.1];\draw(12.5,7) circle [radius=0.1];\draw(12.5,8) circle [radius=0.1];
	\draw(13,5) circle [radius=0.1];\draw(13,6) circle [radius=0.1];\draw(13,7) circle [radius=0.1];\draw(13,8) circle [radius=0.1];
	\draw(13.5,5) circle [radius=0.1];\draw(13.5,6) circle [radius=0.1];\draw(13.5,7) circle [radius=0.1];\draw(13.5,8) circle [radius=0.1];
	\draw(14,5) circle [radius=0.1];\draw(14,6) circle [radius=0.1];\draw(14,7) circle [radius=0.1];\draw(14,8) circle [radius=0.1];
	\draw(14.5,5) circle [radius=0.1];\draw[fill](14.5,6) circle [radius=0.1];\draw(14.5,7) circle [radius=0.1];\draw(14.5,8) circle [radius=0.1];
	\draw(15,5) circle [radius=0.1];\draw(15,6) circle [radius=0.1];\draw(15,7) circle [radius=0.1];\draw(15,8) circle [radius=0.1];
	\draw(15.5,5) circle [radius=0.1];\draw(15.5,6) circle [radius=0.1];\draw(15.5,7) circle [radius=0.1];\draw(15.5,8) circle [radius=0.1];
	\draw(16,5) circle [radius=0.1];\draw(16,6) circle [radius=0.1];\draw(16,7) circle [radius=0.1];\draw(16,8) circle [radius=0.1];
	\draw(16.5,5) circle [radius=0.1];\draw(16.5,6) circle [radius=0.1];\draw(16.5,7) circle [radius=0.1];\draw(16.5,8) circle [radius=0.1];
	\draw(17,5) circle [radius=0.1];\draw(17,6) circle [radius=0.1];\draw(17,7) circle [radius=0.1];\draw(17,8) circle [radius=0.1];
	\draw(17.5,5) circle [radius=0.1];\draw(17.5,6) circle [radius=0.1];\draw(17.5,7) circle [radius=0.1];\draw(17.5,8) circle [radius=0.1];
	\draw(18,5) circle [radius=0.1];\draw(18,6) circle [radius=0.1];\draw(18,7) circle [radius=0.1];\draw(18,8) circle [radius=0.1];
	\draw(18.5,5) circle [radius=0.1];\draw[fill](18.5,6) circle [radius=0.1];\draw(18.5,7) circle [radius=0.1];\draw(18.5,8) circle [radius=0.1];
	\draw(19,5) circle [radius=0.1];\draw(19,6) circle [radius=0.1];\draw(19,7) circle [radius=0.1];\draw(19,8) circle [radius=0.1];
	\draw(19.5,5) circle [radius=0.1];\draw(19.5,6) circle [radius=0.1];\draw(19.5,7) circle [radius=0.1];\draw(19.5,8) circle [radius=0.1];
	\draw(20,5) circle [radius=0.1];\draw(20,6) circle [radius=0.1];\draw(20,7) circle [radius=0.1];\draw(20,8) circle [radius=0.1];
	\draw [thick] (16,4.5) node[below]{$2^{nd}$ layer} ;
	
	\draw [thick] (10,-1.5) node[below]{Figure 6. Quasi-perfect $3$-error correcting code in $ C_4\square C_{16}\square C_4$} ;
	\draw [thick] (10,-2) node[below]{using 3-perfect code in $ C_4\times C_{16}\times C_2$} ;
	\draw [thick] (10,-2.5) node[below]{(filled circles are codewords)} ;
	\draw(6,0)--(6,3);\draw(6.5,0)--(6.5,3);\draw(7,0)--(7,3);\draw(7.5,0)--(7.5,3);\draw(8,0)--(8,3);\draw(8.5,0)--(8.5,3);\draw(9,0)--(9,3);\draw(9.5,0)--(9.5,3);\draw(10,0)--(10,3);\draw(10.5,0)--(10.5,3);\draw(11,0)--(11,3);\draw(11.5,0)--(11.5,3);\draw(12,0)--(12,3);\draw(12.5,0)--(12.5,3);\draw(13,0)--(13,3);\draw(13.5,0)--(13.5,3);
	\draw(6,0)--(13.5,0);\draw(6,1)--(13.5,1);\draw(6,2)--(13.5,2);\draw(6,3)--(13.5,3);
	
	\draw(6,10)--(6,13);\draw(6.5,10)--(6.5,13);\draw(7,10)--(7,13);\draw(7.5,10)--(7.5,13);\draw(8,10)--(8,13);\draw(8.5,10)--(8.5,13);\draw(9,10)--(9,13);\draw(9.5,10)--(9.5,13);\draw(10,10)--(10,13);\draw(10.5,10)--(10.5,13);\draw(11,10)--(11,13);\draw(11.5,10)--(11.5,13);\draw(12,10)--(12,13);\draw(12.5,10)--(12.5,13);\draw(13,10)--(13,13);\draw(13.5,10)--(13.5,13);
	\draw(6,10)--(13.5,10);\draw(6,11)--(13.5,11);\draw(6,12)--(13.5,12);\draw(6,13)--(13.5,13);
	
	\draw(-.5,5)--(-.5,8);\draw(0,5)--(0,8);\draw(0.5,5)--(0.5,8);\draw(1,5)--(1,8);\draw(1.5,5)--(1.5,8);\draw(2,5)--(2,8);\draw(2.5,5)--(2.5,8);\draw(3,5)--(3,8);\draw(3.5,5)--(3.5,8);\draw(4,5)--(4,8);\draw(4.5,5)--(4.5,8);\draw(5,5)--(5,8);\draw(5.5,5)--(5.5,8);\draw(6,5)--(6,8);\draw(6.5,5)--(6.5,8);\draw(7,5)--(7,8);
	\draw(-.5,5)--(7,5);\draw(-.5,6)--(7,6);\draw(-.5,7)--(7,7);\draw(-.5,8)--(7,8);
	
	\draw(12.5,5)--(12.5,8);\draw(13,5)--(13,8);\draw(13.5,5)--(13.5,8);\draw(14,5)--(14,8);\draw(14.5,5)--(14.5,8);\draw(15,5)--(15,8);\draw(15.5,5)--(15.5,8);\draw(16,5)--(16,8);\draw(16.5,5)--(16.5,8);\draw(17,5)--(17,8);\draw(17.5,5)--(17.5,8);\draw(18,5)--(18,8);\draw(18.5,5)--(18.5,8);\draw(19,5)--(19,8);\draw(19.5,5)--(19.5,8);\draw(20,5)--(20,8);
	\draw(12.5,5)--(20,5);\draw(12.5,6)--(20,6);\draw(12.5,7)--(20,7);\draw(12.5,8)--(20,8);
\end{tikzpicture}
	
	\begin{thm}\label{lm7}
		Let \( m, n \geq 2 \), and let \( e \geq 1 \). If there exists a perfect \( e \)-error-correcting code in \( C_m \square C_n \square C_k \) for \( k = 1, 2 \), then a quasi-perfect \( e \)-error-correcting code exists in each of the following Cartesian products:
		\[
		C_m \square C_n \square C_i, \quad 
		C_{m+1} \square C_n \square C_i, \quad 
		C_m \square C_{n+1} \square C_i, \quad 
		C_{m+1} \square C_{n+1} \square C_i
		\]
		for \( i = 1, 2, 3, 4 \).
	\end{thm}
	
	\begin{proof}
		We present the proof for the case \( k = 2 \); the case \( k = 1 \) follows analogously (see Figures~4 and~5).
		
		Assume that \( D = D_1 \cup D_2 \) is a perfect \( e \)-error-correcting code in \( C_m \square C_n \square C_2 \). Then:
		\begin{itemize}
			\item In \( C_m \square C_n \square C_2 \), there are two layers of \( C_m \square C_n \), denoted by \( (C_m \square C_n)_0 \) and \( (C_m \square C_n)_1 \). The sets \( D_1 \subset (C_m \square C_n)_0 \) and \( D_2 \subset (C_m \square C_n)_1 \) form the codewords of \( D \), and the entire graph is covered by disjoint balls \( B_1, \dots, B_\ell \) of radius \( e \) centered at the codewords.
			
			\item In the extended graphs \( C_{m+1} \square C_n \square C_2 \), \( C_m \square C_{n+1} \square C_2 \), and \( C_{m+1} \square C_{n+1} \square C_2 \), we are effectively adding either one row, one column, or both to each layer. The newly added vertices are all at most distance 1 from some vertex in the original graph \( C_m \square C_n \square C_2 \). Hence, by extending each ball \( B_i \) to radius \( e+1 \), we obtain a covering of the extended graphs. Therefore, \( D \) becomes a quasi-perfect \( e \)-error-correcting code in each of these graphs.
			
			\item Now, consider the graph \( C_m \square C_n \square C_3 \), which has three layers \( (C_m \square C_n)_0 \), \( (C_m \square C_n)_1 \), and \( (C_m \square C_n)_2 \). Assume \( D_1 \subset (C_m \square C_n)_0 \), \( D_2 \subset (C_m \square C_n)_1 \). The code \( D = D_1 \cup D_2 \) has minimum distance \( 2e + 1 \). All vertices in layers 0 and 1 are covered within radius \( e \), as \( D \) is a perfect code in \( C_m \square C_n \square C_2 \).
			
			Vertices in layer 2 of the form \( (x,y,2) \), where \( (x,y,0) \in S_e(z) \) for some \( z \in D_1 \), or \( (x,y,1) \in S_e(w) \) for some \( w \in D_2 \), are at distance \( e+1 \) from \( D \). All remaining vertices in layer 2 are within distance \( e \) from some codeword. Thus, the covering radius is \( e+1 \), and \( D \) is quasi-perfect in \( C_m \square C_n \square C_3 \), and similarly in the extended graphs with one additional row and/or column.
			
			\item Next, consider \( C_m \square C_n \square C_4 \), which has four layers \( (C_m \square C_n)_i \) for \( i = 0,1,2,3 \). Let \( D_1 \subset (C_m \square C_n)_0 \) and \( D_2 \subset (C_m \square C_n)_2 \), and define \( D = D_1 \cup D_2 \). Then:
			\begin{itemize}
				\item The minimum distance between any two codewords is \( 2e + 1 \).
				\item Vertices \( (x,y,0) \), where \( (x,y,2) \in S_e(z) \) for some \( z \in D_2 \), and vice versa, are at distance \( e+1 \) from \( D \).
				\item Vertices in layers 1 and 3 that lie within a radius \( e \) neighborhood of any \( (x,y,0) \in D_1 \) or \( (x,y,2) \in D_2 \) are at distance at most \( e+1 \) from \( D \).
				\item All other vertices are within distance \( e \) of some codeword.
			\end{itemize}
			Hence, the covering radius is \( e+1 \), and \( D \) is quasi-perfect in \( C_m \square C_n \square C_4 \), and by extension, in \( C_{m+1} \square C_n \square C_4 \), \( C_m \square C_{n+1} \square C_4 \), and \( C_{m+1} \square C_{n+1} \square C_4 \).
		\end{itemize}
		This completes the proof.
	\end{proof}
	
	\section{\bf Quasi-Perfect Codes in $C_n \square C_n \square C_l$ from Quasi-Perfect Codes in $C_n \square C_n$}
	
	In this section, we construct quasi-perfect $e$-error-correcting codes in the Cartesian product \( C_n \square C_n \square C_l \) for \( e \leq 2 \), by leveraging quasi-perfect $e$-error-correcting codes in \( C_n \square C_n \). 
	
	Explicit constructions for quasi-perfect 1-error-correcting codes in \( C_n \square C_n \square C_n \) for \( n = 3, 4 \) are illustrated in Figures 7 and 8, respectively.
	
	\begin{nt} \label{lm8}\cite{albd} The following results were proved:
		\begin{itemize}
			\item Quasi-perfect 1-error correcting codes exist in \( C_n \square C_n \square C_n \) for \( 8 \leq n \leq 12 \). For \( n = 8, 9 \), the code is \( \{(i, 2i) : 0 \leq i < n\} \); and for \( 10 \leq n \leq 12 \), the code is \( \{(2i, 3i) : 0 \leq i < n\} \).
			
			\item Quasi-perfect 2-error correcting codes exist in \( C_n \square C_n \square C_n \) for \( 14 \leq n \leq 24 \). For \( 14 \leq n \leq 19 \), the code is \( \{(2i, 3i) : 0 \leq i < n\} \); and for \( 20 \leq n \leq 24 \), the code is \( \{(3i, 4i) : 0 \leq i < n\} \).
		\end{itemize}
	\end{nt}
	
	Using these constructions, we now build quasi-perfect codes in \( C_n \square C_n \square C_l \). Let \( D_0 \) denote a quasi-perfect code in \( C_n \square C_n \), as given in Note~\ref{lm8}.
	
	$$\begin{array}{cc}
		\begin{tikzpicture}[scale=0.5]
			\draw (0,0) circle [radius=0.1];\draw  (0,1) circle [radius=0.1];\draw [fill] (0,2) circle [radius=0.1];
			\draw (0,0) -- (0,1) --  (0,2) ;\draw (0,0) -- (1,0) --  (2,0) ;
			\draw  (1,0) circle [radius=0.1];\draw   (1,1)circle [radius=0.1];\draw  (1,2)circle [radius=0.1];
			\draw  (1,0)--   (1,1)--  (1,2);\draw  (0,1)--   (1,1)--  (2,1);
			\draw   (2,0) circle [radius=0.1];\draw  (2,1) circle [radius=0.1];\draw  (2,2) circle [radius=0.1];
			\draw   (2,0)-- (2,1)-- (2,2);\draw   (0,2)-- (1,2)-- (2,2);

			\draw (-2,-4) circle [radius=0.1];\draw  (-2,-3) circle [radius=0.1];\draw  (-2,-2) circle [radius=0.1];
			\draw (-2,-4) --  (-2,-3) -- (-2,-2) ;\draw (-4,-2) --  (-3,-2) -- (-2,-2) ;
			\draw  (-3,-4) circle [radius=0.1];\draw   [fill](-3,-3)circle [radius=0.1];\draw  (-3,-2)circle [radius=0.1];
			\draw  (-3,-4) --(-3,-3)-- (-3,-2);\draw  (-4,-3) --(-3,-3)-- (-2,-3);
			\draw   (-4,-4) circle [radius=0.1];\draw  (-4,-3) circle [radius=0.1];\draw  (-4,-2) circle [radius=0.1];
			\draw   (-4,-4) -- (-4,-3) --  (-4,-2) ;\draw   (-4,-4) -- (-3,-4) --  (-2,-4) ;

			\draw (4,-4) circle [radius=0.1];\draw  (4,-3) circle [radius=0.1];\draw  (4,-2) circle [radius=0.1];
			\draw (4,-4) -- (4,-3) --  (4,-2) ;\draw (4,-4) -- (5,-4) --  (6,-4) ;
			\draw  (5,-4) circle [radius=0.1];\draw   (5,-3)circle [radius=0.1];\draw  (5,-2)circle [radius=0.1];
			\draw  (5,-4)--  (5,-3)--  (5,-2);\draw  (4,-3)--  (5,-3)--  (6,-3);
			\draw   [fill](6,-4) circle [radius=0.1];\draw  (6,-3) circle [radius=0.1];\draw  (6,-2) circle [radius=0.1];
			\draw   (6,-4) --  (6,-3) --  (6,-2);\draw   (4,-2) --  (5,-2) --  (6,-2);
			
			\draw [<->] (-2,-1) arc (180:90:40pt);
			\draw [<->] (-0.5,-3) --(2.5,-3);
			\draw [<->] (4,-1) arc (0:90:40pt);

			\draw [thick] (1,-5.5) node[below]{Figure 7. Quasi-perfect $1$-error correcting } ;
			\draw [thick] (1.7,-6.3) node[below]{ code in $ C_3\square C_3\square C_3$} ;
			\draw [thick] (1.5,-7) node[below]{(filled circles are codewords)} ;
			
		\end{tikzpicture} 
		
		& \begin{tikzpicture}[scale=0.5]
			\draw (0,0) circle [radius=0.1];\draw  (0,1) circle [radius=0.1];\draw  (0,2) circle [radius=0.1];\draw [fill](0,3) circle [radius=0.1];
			\draw (0,0) -- (0,3);\draw (0,0) --  (3,0) ;
			\draw  (1,0) circle [radius=0.1];\draw   (1,1)circle [radius=0.1];\draw  (1,2)circle [radius=0.1];\draw  (1,3) circle [radius=0.1];
			\draw  (1,0)--   (1,3);\draw  (0,1)--   (3,1);
			\draw   (2,0) circle [radius=0.1];\draw  (2,1) circle [radius=0.1];\draw  [fill](2,2) circle [radius=0.1];\draw  (2,3) circle [radius=0.1];
			\draw   (2,0)--  (2,3);\draw   (0,2)-- (3,2);
			\draw   (3,0) circle [radius=0.1];\draw  (3,1) circle [radius=0.1];\draw  (3,2) circle [radius=0.1];\draw  (3,3) circle [radius=0.1];
			\draw   (3,0)--  (3,3);\draw   (0,3)-- (3,3);
			
			\draw (-1,-4) circle [radius=0.1];\draw  (-1,-3) circle [radius=0.1];\draw  (-1,-2) circle [radius=0.1];\draw  (-1,-1) circle [radius=0.1];
			\draw (-1,-4) --  (-1,-1) ;\draw (-1,-4)  -- (-4,-4) ;
			\draw  [fill](-2,-4) circle [radius=0.1];\draw   (-2,-3)circle [radius=0.1];\draw  (-2,-2)circle [radius=0.1];\draw  (-2,-1) circle [radius=0.1];
			\draw  (-2,-4) -- (-2,-1);\draw  (-1,-3) -- (-4,-3);
			\draw   (-3,-4) circle [radius=0.1];\draw  (-3,-3) circle [radius=0.1];\draw  (-3,-2) circle [radius=0.1];\draw  (-3,-1) circle [radius=0.1];
			\draw   (-3,-4)  --  (-3,-1) ;\draw   (-1,-2) --  (-4,-2) ;
			\draw   (-4,-4) circle [radius=0.1];\draw [fill] (-4,-3) circle [radius=0.1];\draw  (-4,-2) circle [radius=0.1];\draw  (-4,-1) circle [radius=0.1];
			\draw   (-4,-4)  --  (-4,-1) ;\draw   (-1,-1) --   (-4,-1) ;
			
			\draw (4,-4) circle [radius=0.1];\draw  (4,-3) circle [radius=0.1];\draw  (4,-2) circle [radius=0.1];\draw  (4,-1) circle [radius=0.1];
			\draw (4,-4)  --  (4,-1) ;\draw (4,-4)  --  (7,-4) ;
			\draw  (5,-4) circle [radius=0.1];\draw   [fill](5,-3)circle [radius=0.1];\draw  (5,-2)circle [radius=0.1];\draw  (5,-1) circle [radius=0.1];
			\draw  (5,-4)--  (5,-1);\draw  (4,-3)--  (7,-3);
			\draw   (6,-4) circle [radius=0.1];\draw  (6,-3) circle [radius=0.1];\draw  (6,-2) circle [radius=0.1];\draw  (6,-1) circle [radius=0.1];
			\draw   (6,-4)  --  (6,-1);\draw   (4,-2)  --  (7,-2);
			\draw   [fill](7,-4) circle [radius=0.1];\draw  (7,-3) circle [radius=0.1];\draw  (7,-2) circle [radius=0.1];\draw  (7,-1) circle [radius=0.1];
			\draw   (7,-4)  --  (7,-1);\draw   (4,-1)  --  (7,-1);

			\draw (0,-8) circle [radius=0.1];\draw  (0,-7) circle [radius=0.1];\draw  (0,-6) circle [radius=0.1];\draw  (0,-5) circle [radius=0.1];
			\draw (0,-8) --   (0,-5) ;\draw (0,-8) --   (3,-8) ;
			\draw (1,-8) circle [radius=0.1];\draw  (1,-7) circle [radius=0.1];\draw  (1,-6) circle [radius=0.1];\draw [fill] (1,-5) circle [radius=0.1];
			\draw (1,-8) --   (1,-5) ;\draw (0,-7) --   (3,-7) ;
			\draw (2,-8) circle [radius=0.1];\draw  (2,-7) circle [radius=0.1];\draw  (2,-6) circle [radius=0.1];\draw  (2,-5) circle [radius=0.1];
			\draw (2,-8) --   (2,-5) ;\draw (0,-6) --   (3,-6) ;
			\draw (3,-8) circle [radius=0.1];\draw  (3,-7) circle [radius=0.1];\draw  [fill](3,-6) circle [radius=0.1];\draw  (3,-5) circle [radius=0.1];
			\draw (3,-8) --   (3,-5) ;\draw (0,-5) --   (3,-5) ;
			
			\draw [<->] (-2,0) arc (180:90:20pt);
			\draw [<->] (5,0) arc (0:90:20pt);
			\draw [<->] (-2,-5) arc (180:270:20pt);
			\draw [<->] (5,-5) arc (360:270:20pt);
			
			\draw [thick] (1,-9.5) node[below]{Figure 8. Quasi-perfect $1$-error correcting } ;
			\draw [thick] (1.7,-10.3) node[below]{code in $ C_4\square C_4\square C_4$} ;
			\draw [thick] (1.5,-11) node[below]{(filled circles are codewords)} ;
			
		\end{tikzpicture}
	\end{array}$$
	
	\begin{thm}\label{lm9}
		There exists a perfect 1-error correcting code in the Cartesian product \( C_6 \square C_6 \square C_2 \).
	\end{thm}
	
	\begin{proof}
		Let \( D_0 \) be a quasi-perfect 1-error correcting code in \( C_6 \square C_6 \), and define \( D_1 = (0,3) + D_0 \). Then 
		\[
		D = \bigcup_{i=0}^{1} D_i \oplus \{i\}
		\]
		is a perfect 1-error correcting code in \( C_6 \square C_6 \square C_2 \). The minimum distance of \( D \) is 3, and the covering radius is 1.
	\end{proof}
	
	\cite{go} proved the existence of a perfect 2-error correcting code in \( C_7 \square C_7 \square C_7 \).
	
	\begin{thm}\label{lm10}
		For any integer \( k \geq 1 \), there exists a quasi-perfect 1-error correcting code in the Cartesian product of:
		\begin{enumerate}
			\item \( C_6 \square C_6 \square C_{3k} \),
			\item \( C_n \square C_n \square C_{3k} \) for \( 8 \leq n \leq 12 \),
			\item \( C_n \square C_n \square C_n \) for \( 8 \leq n \leq 12 \).
		\end{enumerate}
	\end{thm}
	
	\begin{proof}
		\begin{enumerate}
			\item \textbf{For \( C_6 \square C_6 \square C_{3k} \):}  
			Define:
			\[
			D_0 = \text{a quasi-perfect 1-error correcting code in } C_6 \square C_6,
			\]
			\[
			D_1 = \{(0,3), (2,1), (4,5)\}, \quad D_2 = \{(1,5), (3,3), (5,1)\}.
			\]
			Then the code
			\[
			D = \bigcup_{i=0}^{2} D_i \oplus \{i\}
			\]
			is a quasi-perfect 1-error correcting code in \( C_6 \square C_6 \square C_3 \). The vertices of \( D_1 \oplus \{1\} \) and \( D_2 \oplus \{2\} \) are at distance 3 from \( D_0 \oplus \{0\} \), ensuring the minimum distance is 3 and covering radius is 2.
			
			Extending this to \( C_6 \square C_6 \square C_{3k} \), define:
			\[
			D = \bigcup_{i=0}^{k-1} \left( D_0 \oplus \{3i\} \cup D_1 \oplus \{3i+1\} \cup D_2 \oplus \{3i+2\} \right).
			\]
			This forms a quasi-perfect 1-error correcting code as above.
			
			\item \textbf{For \( C_n \square C_n \square C_{3k} \):}  
			Let \( D_0 \) be a quasi-perfect 1-error correcting code in \( C_n \square C_n \). Define:
			\[
			D_1 = (0,3) + D_0, \quad D_2 = (0,n-3) + D_0.
			\]
			Then
			\[
			D = \bigcup_{i=0}^{2} D_i \oplus \{i\}
			\]
			is a quasi-perfect 1-error correcting code in \( C_n \square C_n \square C_3 \), with minimum distance 3 and covering radius 2.
			
			Generalizing to \( C_n \square C_n \square C_{3k} \), define:
			\[
			D = \bigcup_{i=0}^{k-1} \left( D_0 \oplus \{3i\} \cup D_1 \oplus \{3i+1\} \cup D_2 \oplus \{3i+2\} \right).
			\]
			
			\item \textbf{For \( C_n \square C_n \square C_n \):}  
			Define:
			\[
			D_i = (0,3i) + D_0, \quad \text{for } 1 \leq i \leq n - 1.
			\]
			Then the code
			\[
			D = \bigcup_{i=0}^{n-1} D_i
			\]
			is a quasi-perfect 1-error correcting code in \( C_n \square C_n \square C_n \), with the same properties.
		\end{enumerate}
	\end{proof}

	\begin{thm}\label{lm11}
		For any integer \( k \geq 1 \), there exists a quasi-perfect \( 2 \)-error correcting code in the Cartesian product:
		\begin{enumerate}
			\item \( C_{14} \square C_{14} \square C_{4k} \),
			\item \( C_n \square C_n \square C_{6k} \) for \( 14 \leq n \leq 19 \).
		\end{enumerate}
	\end{thm}
	
	\begin{proof}
		\begin{enumerate}
			\item \textbf{For \( C_{14} \square C_{14} \square C_{4k} \):}
			
			We first construct a quasi-perfect 2-error correcting code in \( C_{14} \square C_{14} \square C_4 \).
			
			\begin{itemize}
				\item Let \( D_0 \) be a quasi-perfect 2-error correcting code in \( C_{14} \square C_{14} \). By direct computation, we observe that there are exactly 14 vertices in \( C_{14} \square C_{14} \) that lie at distance 3 from \( D_0 \). Define \( D_1 = (1, n-2) + D_0 \); this shifts the original code \( D_0 \) such that these previously uncovered vertices are now covered.
				
				\item Define 
				\[
				D = (D_0 \oplus \{0\}) \cup (D_1 \oplus \{2\}).
				\]
				Then \( D \) is a quasi-perfect 2-error correcting code in \( C_{14} \square C_{14} \square C_4 \). Since the vertices in \( D_1 \oplus \{2\} \) are at distance 5 from those in \( D_0 \oplus \{0\} \), the minimum distance of \( D \) is 5.
				
				\item All vertices in layers \( C_{14} \square C_{14} \times \{0,2\} \) are covered by radius-2 balls centered at codewords in \( D \). The remaining layers \( C_{14} \square C_{14} \times \{1,3\} \) are at distance 3 from the code, so the covering radius is 3.
				
				\item To extend this to \( C_{14} \square C_{14} \square C_{4k} \), define
				\[
				D = \bigcup_{i=0}^{k-1} \left( D_0 \oplus \{4i\} \cup D_1 \oplus \{4i + 2\} \right).
				\]
				By periodic repetition of the code blocks, this gives a quasi-perfect 2-error correcting code in \( C_{14} \square C_{14} \square C_{4k} \).
			\end{itemize}
			
			\item \textbf{For \( C_n \square C_n \square C_{6k} \), where \( 14 \leq n \leq 19 \):}
			
			\begin{itemize}
				\item Let \( D_0 \) be a quasi-perfect 2-error correcting code in \( C_n \square C_n \), as given in Note~\ref{lm8}. Define two additional code sets:
				\[
				D_1 = (1,5) + D_0, \quad D_2 = (3,1) + D_0.
				\]
				\item Define
				\[
				D = \bigcup_{i=0}^{2} \left( D_i \oplus \{2i\} \right).
				\]
				Using a computer search, it was verified that this forms a quasi-perfect 2-error correcting code in \( C_n \square C_n \square C_6 \). The minimum distance is 5, and the covering radius is 3.
				
				\item To generalize to \( C_n \square C_n \square C_{6k} \), define:
				\[
				D = \bigcup_{i=0}^{k-1} \left( D_0 \oplus \{6i\} \cup D_1 \oplus \{6i+2\} \cup D_2 \oplus \{6i+4\} \right).
				\]
				By construction, this yields a quasi-perfect 2-error correcting code in \( C_n \square C_n \square C_{6k} \).
			\end{itemize}
		\end{enumerate}
	\end{proof}
	
	\section{\bf Quasi-perfect Codes in $P_m \square P_n$ and $P_m \square P_n \square P_l$ }
	
	In this section, we study quasi-perfect codes in the Cartesian product of two and three paths, namely \( P_m \square P_n \) and \( P_m \square P_n \square P_l \).
	
	It is easy to observe that for every \( n \geq 2 \), the set
	\[
	D = \{(0,0), (n-1,n-1)\}
	\]
	forms an \((n-2)\)-quasi-perfect code in \( P_n \square P_n \).
	
	\begin{thm}\label{lm12}
		The Cartesian product \( P_m \square P_n \), where \( 2 \leq m \leq n \), admits an \( e \)-quasi-perfect code for \( e \geq 1 \) if one of the following holds:
		\begin{itemize}
			\item \( m = n = 2e + 3 \),
			\item \( m = e + 1 \) and \( n = e + 3 \).
		\end{itemize}
	\end{thm}
	
	\begin{proof}
		We consider each case separately.
		
		\begin{itemize}
			\item \textbf{Case 1:} \( m = n = 2e + 3 \).\\
			Define the code
			\[
			D = \{(1,1), (e+2, e+2), (n,n), (1,n), (n,1)\}.
			\]
			The minimum pairwise distance among the codewords is \( 2e + 2 = n - 1 \). All vertices in \( P_n \square P_n \), except those in the sphere \( S_{e+1}((e+2, e+2)) \), are covered by balls of radius \( e \) centered at codewords in \( D \). The uncovered vertices in \( S_{e+1}((e+2, e+2)) \) are at distance \( e + 1 \) from the nearest codeword. Thus, the covering radius is \( e + 1 \), and \( D \) is an \( e \)-quasi-perfect code.
			
			\item \textbf{Case 2:} \( m = e + 1, ~ n = e + 3 \).\\
			Define the code
			\[
			D = \{(1,1), (m, n-1)\}.
			\]
			The minimum distance between the two codewords is \( 2e + 1 \). Every vertex in \( P_m \square P_n \), except \( (1,n) \), lies within distance \( e \) of some codeword. The vertex \( (1,n) \) is at distance \( e + 1 \) from both codewords. Therefore, the covering radius is \( e + 1 \), and \( D \) is an \( e \)-quasi-perfect code.
		\end{itemize}
	\end{proof}
	
	\begin{obs}
		For all \( n \geq 2 \), the Cartesian product \( P_n \square P_n \square P_2 \) admits a perfect \((n-1)\)-error correcting code. One such code is:
		\[
		D = \{(0,0,0), (n-1, n-1, 1)\}.
		\]
	\end{obs}
	\begin{tikzpicture}[scale=0.4]
		
		\draw  [fill](1,14) circle [radius=0.4];\draw [fill](17,-2) circle [radius=0.4];
		
		\draw[fill](-.2,-.2)rectangle(.2,.2);\draw[fill](2.8,-3.2)rectangle(3.2,-2.8);\draw[fill](2.8,1.8)rectangle(3.2,2.2);
		\draw[fill](4.8,-.2)rectangle(5.2,.2);\draw[fill](4.8,4.8)rectangle(5.2,5.2);\draw[fill](7.8,1.8)rectangle(8.2,2.2);
		\draw[fill](7.8,6.8)rectangle(8.2,7.2);\draw[fill](9.8,4.8)rectangle(10.2,5.2);\draw[fill](9.8,9.8)rectangle(10.2,10.2);
		\draw[fill](12.8,6.8)rectangle(13.2,7.2);\draw[fill](12.8,11.8)rectangle(13.2,12.2);\draw[fill](14.8,9.8)rectangle(15.2,10.2);\draw[fill](14.8,14.8)rectangle(15.2,15.2);
		\draw[fill](17.8,11.8)rectangle(18.2,12.2);

		\draw(0,0)--(0,15);\draw(0,0)--(15,0);
		\draw(1,-1)--(1,14);\draw(1,-1)--(16,-1);
		\draw(2,-2)--(2,13);\draw(2,-2)--(17,-2);
		\draw(3,-3)--(3,12);\draw(3,-3)--(18,-3);
		\draw (5,0)--(5,15);\draw (6,-1) --(6,14);\draw (7,-2)--(7,13);\draw  (8,-3)--(8,12);
		\draw (10,0)--(10,15);\draw (11,-1)--(11,14);\draw (12,-2)--(12,13);\draw  (13,-3) -- (13,12) ;
		\draw (15,0)--(15,15);\draw (16,-1)--(16,14);\draw (17,-2)--(17,13);\draw  (18,-3)--(18,12);

		\draw(3,2)--(18,2);\draw(3,7)--(18,7);\draw(3,12)--(18,12);
		\draw(2,3)--(17,3);\draw(2,8)--(17,8);\draw(2,13)--(17,13);
		\draw(1,4)--(16,4);\draw(1,9)--(16,9);\draw(1,14)--(16,14);
		\draw(0,5)--(15,5);\draw(0,10)--(15,10);\draw(0,15)--(15,15);
		
		\draw(0,0)--(3,-3);\draw(5,0)--(8,-3);\draw(10,0)--(13,-3);\draw(15,0)--(18,-3);
		\draw(0,5)--(3,2);\draw(5,5)--(8,2);\draw(10,5)--(13,2);\draw(15,5)--(18,2);
		\draw(0,10)--(3,7);\draw(5,10)--(8,7);\draw(10,10)--(13,7);\draw(15,10)--(18,7);
		\draw(0,15)--(3,12);\draw(5,15)--(8,12);\draw(10,15)--(13,12);\draw(15,15)--(18,12);

		
		\draw [thick] (7,-4.5) node[below]{Figure 9. Quasi perfect $3$-error correcting code in $ P_4\square P_4\square P_4$} ;
		\draw [thick] (7,-5.5) node[below]{(filled circles are codewords)} ;
		\draw [thick] (7,-6.5) node[below]{(filled squares are vertices at a distance $3$ from at least one codeword)} ;
	\end{tikzpicture}

	\begin{tikzpicture}[scale=0.4]
		
		\draw[fill] (0,15) circle [radius=0.4];\draw[fill] (17,-2) circle [radius=0.4];
		\draw[fill](.8,-1.2)rectangle(1.2,-.7);\draw[fill](1.8,2.8)rectangle(2.2,3.2);\draw[fill](4.8,-.2)rectangle(5.2,.2);
		\draw[fill](5.8,3.8)rectangle(6.2,4.2);\draw[fill](6.8,7.8)rectangle(7.2,8.2);\draw[fill](9.8,4.8)rectangle(10.2,5.2);
		\draw[fill](10.8,8.8)rectangle(11.2,9.2);\draw[fill](11.8,12.8)rectangle(12.2,13.2);
		\draw[fill](14.8,9.8)rectangle(15.2,10.2);\draw[fill](15.8,13.8)rectangle(16.2,14.2);

		\draw(0,0)--(0,15);\draw(0,0)--(15,0);\draw(1,-1)--(1,14);\draw(1,-1)--(16,-1);\draw(2,-2)--(2,13);
		\draw(2,-2)--(17,-2);
		\draw (5,0)--(5,15);\draw (6,-1) --(6,14);\draw (7,-2)--(7,13);
		\draw (10,0)--(10,15);\draw (11,-1)--(11,14);\draw (12,-2)--(12,13);
		\draw (15,0)--(15,15);\draw (16,-1)--(16,14);\draw (17,-2)--(17,13);

		\draw(2,3)--(17,3);\draw(2,8)--(17,8);\draw(2,13)--(17,13);
		\draw(1,4)--(16,4);\draw(1,9)--(16,9);\draw(1,14)--(16,14);
		\draw(0,5)--(15,5);\draw(0,10)--(15,10);\draw(0,15)--(15,15);
		
		\draw(0,0)--(2,-2);\draw(5,0)--(7,-2);\draw(10,0)--(12,-2);\draw(15,0)--(17,-2);
		\draw(0,5)--(2,3);\draw(5,5)--(7,3);\draw(10,5)--(12,3);\draw(15,5)--(17,3);
		\draw(0,10)--(2,8);\draw(5,10)--(7,8);\draw(10,10)--(12,8);\draw(15,10)--(17,8);
		\draw(0,15)--(2,13);\draw(5,15)--(7,13);\draw(10,15)--(12,13);\draw(15,15)--(17,13);

		
		\draw [thick] (7,-3.5) node[below]{Figure 10. Quasi perfect $3$-error correcting code in $ P_4\square P_4\square P_3$} ;
		\draw [thick] (7,-4.5) node[below]{(filled circles are codewords)} ;
		\draw [thick] (7,-5.5) node[below]{(filled squares are vertices at a distance $3$ from at least one codeword)} ;
	\end{tikzpicture}

	\begin{thm}\label{lm13}
		For all \( n \geq 2 \), there exists an \((n-1)\)-quasi-perfect code in the Cartesian products \( P_n \square P_n \square P_3 \) and \( P_n \square P_n \square P_4 \).
	\end{thm}
	
	\begin{proof}
		Define \( D = \{(0,0,1), (n-1,n-1,2)\} \). Then for all \( n \geq 2 \), \( D \) is an \((n-1)\)-quasi-perfect code in \( P_n \square P_n \square P_4 \) (see Fig.~9). Similarly, define \( D = \{(0,0,0), (n-1,n-1,2)\} \); this forms an \((n-1)\)-quasi-perfect code in \( P_n \square P_n \square P_3 \) (see Fig.~10). In both cases, the minimum distance between codewords is \( 2(n - 1) + 1 \), and all vertices are within distance \( n \) from the nearest codeword. The covering radius is therefore \( n \), and every vertex lies within distance \( n \) or less from some codeword. Thus, \( D \) is an \((n-1)\)-quasi-perfect code in the respective Cartesian products.
	\end{proof}
	
	\vspace{0.5em}
	\noindent \textbf{Concluding Remarks.} \\
	We have shown that quasi-perfect \( e \)-error correcting codes can be constructed in the Cartesian product of a graph \( G \) with a path or cycle, provided that a perfect \( e \)-error correcting code exists in \( G \). For \( m, n \geq 3 \), we explicitly constructed quasi-perfect 2-error correcting codes in \( P_m \square P_n \square P_{6k - 2} \) and \( C_m \square C_n \square C_{6k} \) for all integers \( k \geq 1 \), based on perfect 2-error correcting codes in \( P_m \square P_n \) and \( C_m \square C_n \), respectively.
	
	Additionally, we constructed quasi-perfect codes in \( P_4 \square P_4 \square P_4 \) by using a perfect code in \( P_2 \square P_2 \square P_2 \). Quasi-perfect codes were also developed in \( C_n \square C_n \square C_l \) for \( 3 \leq n \leq 19 \) and suitable values of \( l \), utilizing known quasi-perfect codes in \( C_n \square C_n \).
	
	A natural direction for further research is to determine, for which integers \( n \) and for which graphs \( G_2 \), one can construct quasi-perfect codes in the Cartesian product of \( G_1 \) with \( n \) copies of \( G_2 \), i.e., in \( G_1 \square G_2 \square \cdots \square G_2 \). Additionally, it would be of interest to identify all values of \( m, n, l \) for which \( C_m \square C_n \square C_l \) admits a quasi-perfect code.

	\section*{Acknowledgment}
	The authors would like to express their sincere gratitude to the anonymous referee for their valuable suggestions and insightful comments, which have significantly improved the quality and presentation of this paper.\\
	
	The first author gratefully acknowledges the Department of Science and Technology, New Delhi, India, for awarding the Women Scientist Scheme (DST/WOS-A/PM-14/2021(G)) for research in Basic/Applied Sciences.

	{\centerline{************}}

	\bibliographystyle{amsplain}

\begin{thebibliography}{10}		
		
		\bibitem{al}  B. AlBdaiwi and P. Horak, “Perfect distance-d placements in 3-dimensional tori,” J. Combinat. Math. Combinat. Comput., vol. 43, pp. 159–174, 2002.
		\bibitem{alb} B. AlBdaiwi and B. Bose, “Quasi-perfect Lee distance codes,” IEEE Trans. Inf. Theory, vol. 49, no. 6, pp. 1535–1539, 2003.
		\bibitem{albda} B. AlBdaiwi and M. Livingston, "Perfect distance-d placements in 2-dimensional tori," J. Supercomput., vol. 29, no. 1, pp.45-57, 2004.
		\bibitem{albdai} B. AlBdaiwi and B. Bose, "Quasi-perfect resource placements for two-dimensional toroidal networks," J. Parallel Distrib. Comput., vol. 65, pp. 815-831, 2005.
		\bibitem{albd} B. AlBdaiwi, P. Horak, and L. Milazzo, “Enumerating and decoding perfect linear Lee codes,” Designs, Codes Cryptogr., vol. 52, no. 2, pp. 155–162, 2009.
		\bibitem{ba} T. Baicheva, I. Bouyukliev, S. Dodunekov, and V. Fack., "Binary and Ternary Linear Quasi-Perfect Codes With Small Dimensions.," IEEE Trans. Inform. Theory, 54(9):4335–4339, 2008.
		\bibitem{bi} K. Bibak, B. Kapron, and V. Srinivasan, “The Cayley graphs associated with some quasi-perfect Lee codes are Ramanujan graphs,” IEEE Trans. Inf. Theory, vol. 62, no. 11, pp. 6355–6358, Nov. 2016.
		\bibitem{br} R. Brualdi, V. Pless, and R. Wilson, “Short codes with a given covering radius,” IEEE Trans. Inf. Theory, vol. 35, no. 1, pp. 99–109,	Jan. 1989.
		\bibitem{bu} S. Buzaglo and T. Etzion, “Tilings by (0.5, n)-crosses and perfect codes,” SIAM J. Discrete Math., vol. 27, no. 2, pp. 1067–1081, 2013.
		\bibitem{ca} C. Camarero and C. Martínez, “Quasi-perfect Lee codes of radius 2 and arbitrarily large dimension,” IEEE Trans. Inf. Theory, vol. 62, no. 3, pp. 1183–1192, Mar. 2016.
		\bibitem{da} D. Danev and S. Dodunekov., "A family of ternary quasi-perfect BCH codes.," Des. Codes Cryptogr., 49(1-3):265–271, 2008.
		\bibitem{dav} A. Davydov and A. Drozhzhina-Labinskaya, “Constructions,
		families, and tables of binary linear covering codes,” IEEE Trans. Inf. Theory, vol. 40, no. 4, pp. 1270–1279, Jul. 1994.
		\bibitem{davy} A. Davydov and L. Tombak, “Quasi-perfect linear binary codes with distance 4 and complete caps in projective geometry,” Probl. Inf.	Transm., vol. 25, no. 4, pp. 265–275, 1989.
		\bibitem{et} T. Etzion and G. Greenberg, “Constructions for perfect mixed codes	and other covering codes,” IEEE Trans. Inf. Theory, vol. 39, no. 1, pp. 209–214, Jan. 1993.
		\bibitem{etz} T. Etzion and B. Mounits., "Quasi-perfect codes with small distance," IEEE Trans. Inform. Theory, 51(11):3928–3946, 2005.
		\bibitem{gab} E. Gabidulin, A. Davydov, and L. Tombak, “Linear codes with covering radius 2 and other new covering codes,” IEEE Trans. Inf. Theory, vol. 37, no. 1, pp. 219–224, Jan. 1991.
		\bibitem{go} S. Golomb and L. Welch, “Perfect codes in the Lee metric and the packing of polyominoes,” SIAM J. Appl. Math., vol. 18, no. 2, pp. 302–317, 1970.
		\bibitem{gor} D. Gorenstein, W. Peterson, and N. Zierler, “Two-error correcting Bose-Chaudhuri codes are quasi-perfect,” Inf. Contr., vol. 3, pp.	291–294, 1960.
		\bibitem{ho} P. Horak and B. AlBdaiwi, “Fast decoding of quasi-perfect Lee distance codes,” Designs, Codes Cryptogr., vol. 40, no. 3, pp. 357–367,	2006.
		\bibitem{hor} P. Horak and B. AlBdaiwi, “Non-periodic tilings of $R^n$ by crosses,” Discrete Comput. Geometry, vol. 47, no. 1, pp. 1–16, 2012.
		\bibitem{hora} P. Horak and V. Hromada, “Tiling $R^5$ by crosses,” Discrete Comput. Geometry, vol. 51, no. 2, pp. 269–284, 2014.
		\bibitem{horak} P. Horak and O. Grošek, “A new approach towards the Golomb–Welch conjecture,” Eur. J. Combinat., vol. 38, pp. 12–22, May 2014.
		\bibitem{le} C. Lee, “Some properties of nonbinary error-correcting codes,” IRE Trans. Inf. Theory, vol. IT-4, no. 2, pp. 77–82, 1958.
		\bibitem{li} M. Livingston and Q. Stout, "Perfect dominating sets," Congr. Numer., vol. 79, pp. 187-203, 1990.
		\bibitem{mo} O. Moreno, “Further results on quasi-perfect codes related to the Goppa codes,” Congressus Numerantium, vol. 40, pp. 249–256, 1983.
		\bibitem{qu} C. Queiroz, C. Camarero, C. Martínez, and R. Palazzo, Jr., “Quasi-perfect codes from Cayley graphs over integer rings,” IEEE Trans. Inf. Theory, vol. 59, no. 9, pp. 5905–5916, Sep. 2013.
		\bibitem{st} J. Strapasson, G. Jorge, A. Campello. and S. Costa, "Quasi-perfect codes in the $l_p$ metric," Comp. Appl. Math., vol. 37, pp. 852–866, 2018.
		\bibitem{str} R. Struik, “Covering Codes,” Ph.D. dissertation, Eindhoven Univ.	Technol., Eindhoven, The Netherlands, 1994.	
		\bibitem{sz} S. Szabó, “On mosaics consisting of multidimensional crosses,” Acta 	Math. Acad. Sci. Hung., vol. 38, nos. 1–4, pp. 191–203, 1981.
		\bibitem{wa} T. Wagner, “A search technique for quasi-perfect codes,” Inf. Contr., vol. 9, pp. 94–99, 1966.
		\bibitem{we} D. West, \textit{Introduction to Graph Theory}, $2^{nd}$ ed., Prentice Hall of India, New Delhi, 2002. 
		\bibitem{za} G. Zaitsev, V. Zinovjev, and N. Semakov, “Interrelation of Preparata and Hamming codes and extension of Hamming codes to
		new double-error-correcting codes,” in Proc. 2nd Int. Symp. Information Theory, B. N. Petrov and F. Csáki, Eds. Budapest, Hungary:
		Akadémiai Kiadó, 1973, pp. 257–263.
		\bibitem{zh} T. Zhang and G. Ge, “Perfect and quasi-perfect codes under the $l_p$ metric,” IEEE Trans. Inf. Theory, vol. 63, no. 7, pp. 4325–4331, Jul. 2017.
		
		
		{\centerline{************}}
		
		
		
		
		
		
	\end{thebibliography}

\end{document}